\documentclass[a4paper,12pt]{article}
\usepackage{amsmath}
\usepackage{amssymb}
\usepackage{tabularx}
\usepackage{enumerate}
\usepackage{graphicx,subfigure}
\newtheorem{thm}{Theorem}[section]
\newtheorem{lem}[thm]{Lemma}

\newtheorem{defi}{Definition}[section]
\newtheorem{pppp}{Proof}

\newcommand{\qed}{\hspace{1em}\mbox{\raisebox{0.65ex}{\fbox{}}}}

\numberwithin{equation}{section}

\newcommand{\be}{\begin{equation}}
\newcommand{\ee}{\end{equation}}
\newcommand\bes{\begin{eqnarray}} \newcommand\ees{\end{eqnarray}}
\newcommand{\bess}{\begin{eqnarray*}}
\newcommand{\eess}{\end{eqnarray*}}

\newcommand{\bpf}{{\bf Proof:\ \ }}
\newcommand{\epf}{\mbox{}\hfill $\Box$}

\begin{document}

\bibliographystyle{plain}
\baselineskip = 20pt
\thispagestyle{empty}

\title{{\bf Analysis of a non-autonomous mutualism
model driven by  L$\acute{e}$vy jumps\thanks{The work is supported by Supported in part by a NSFC Grant No. 11171158,
NSF of Jiangsu Education Committee No. 11KJA110001,  PAPD of Jiangsu Higher Education Institutions and the Jiangsu
Collaborative Innovation Center for Climate Change and Project of Graduate Education Innovation of Jiangsu Province No. KYLX\underline{ }0719.}}}
\date{\empty}
\author{Mei Li$^{a, b}$, Hongjun Gao$^{a,}$\footnote{The corresponding author. E-mail: gaohj@njnu.edu.cn} , Binjun Wang$^{a}$ \\
{\small $^a$ Institute of mathematics, Nanjing Normal University,}\\
{\small Nanjing 210023, PR China}\\
{\small $^b$ School of Applied Mathematics, Nanjing University of Finance and Economics,}\\
{\small Nanjing 210023, PR China}\\
{\small Email: limei@njue.edu.cn }}

\maketitle

\begin{quote}
\noindent
{\bf Abstract.} { %\footnotesize
\small This article is concerned with a mutualism ecological model with L$\acute{e}$vy noise. The local existence and
uniqueness of a positive solution are obtained with positive initial value, and the asymptotic  behavior
to the  problem is studied. Moreover, we show that the solution is stochastically bounded and stochastic permanence.
The sufficient conditions for the system to be extinct are given and the conditions for the system
to be persistence in mean are also established.
}

\end{quote}

\section{Introduction}

Mutualism is an important biological interaction in nature. It
occurs when one species provides some benefit in exchange for some
benefit, for example, pollinators and flowering plants,  the
pollinators obtain floral nectar (and in some cases pollen) as a
food resource while the plant obtains non-trophic reproductive
benefits through pollen dispersal and seed production. Another
instance is ants and aphids, in which the ants obtain honeydew
food resources excreted by aphids while the aphids obtain
increased survival by the non-trophic service of ant defense
against natural enemies of the aphids. Lots of authors have
discussed these models \cite{AR, BJK, CC, CY, GB, HD1, HDB, HS, LT,
TNS}. One of the simplest models is the classical Lotka-Volterra
two-species mutualism model as follows:
\begin{eqnarray}
\left\{
\begin{array}{lll}
\dot{x}(t)=x(t)\big(a_1-b_1x(t)+c_1y(t)\big), \\
\dot{y}(t)=y(t)\big(a_2-b_2y(t)+c_2x(t)\big).
\end{array} \right.
\label{a1}
\end{eqnarray}

Among various types mutualistic model, we should specially mention the following model which was proposed by May \cite{MRM} in 1976:
\begin{eqnarray}
\left\{
\begin{array}{ll}
\dot{x}(t)=x(t)\big(r_1-\frac{b_1x(t)}{K_1+y(t)}-\varepsilon_1x(t)\big),& \\
\dot{y}(t)=y(t)\big(r_2-\frac{b_2y(t)}{K_2+x(t)}-\varepsilon_2y(t)\big),&
\end{array} \right.
\label{aode1}
\end{eqnarray}
where $x(t), y(t)$ denote population densities of each species at time t, $r_i, K_i, b_i, \varepsilon_i$ (i=1, 2)
are positive constants, $r_1, r_2$ denote the intrinsic growth rate of species $x(t), y(t)$ respectively, $K_1$ is the capability of species $x(t)$ being short of  $y(t)$, similarly $K_2$ is the capability of species $y(t)$ being short of  $x(t)$. For \eqref{aode1}, there are three trivial equilibrium points
$$ E_1=(0,0),\ \ E_2=(\frac{r_1}{\varepsilon_1+\frac{b_1}{K_1}}, 0),\ \ E_3=(0, \frac{r_2}{\varepsilon_2+\frac{b_2}{K_2}}),$$
and a unique positive interior equilibrium point $E^*=(x^*, y^*)$ satisfying the following equations
\begin{eqnarray}
\left\{
\begin{array}{ll}
r_1-\frac{b_1x(t)}{K_1+y(t)}-\varepsilon_1x(t)=0,& \\
r_2-\frac{b_2y(t)}{K_2+x(t)}-\varepsilon_2y(t)=0,&
\end{array} \right.
\label{aode2}
\end{eqnarray}
where $E^*$ is globally asymptotically stable.

 In addition, population dynamics is inevitably affected by environmental noises,
  May \cite{MR} pointed out the fact that due to environmental fluctuation, the birth rates, carrying capacity,  and
  other parameters involved in the model system exhibit random fluctuation to a greater or lesser extent. Consequently the equilibrium population
   distribution fluctuates randomly around some average values. Therefore lots of authors introduced stochastic perturbation into deterministic models
   to reveal the effect of environmental variability on the population dynamics in mathematical ecology \cite{DS, HWH, JJ, JJS, JSL, LGJM, LW1, LW2, LW3, MSR, TDHS}.
   Li and Gao et al took into account the effect of randomly fluctuating environment in \cite{LG}, where
   they  considered white noise to each equation of the problem \eqref{aode1}.
   Suppose that parameter $r_i$ is stochastically perturbed, with
   $$r_i\rightarrow r_i+\alpha_i\dot{W}_i(t), \  i=1, 2,$$
   where $W_1(t), W_2(t)$ are mutually independent Brownian motion, $\alpha_i, \ i=1,2 $ represent the intensities of the white noise. Then
   the corresponding deterministic model system \eqref{aode1} may be described by the It$\hat{o}$ problems:
\begin{eqnarray}
\left\{
\begin{array}{ll}
dx(t)=x(t)\big(r_1-\frac{b_1x(t)}{K_1+y(t)}-\varepsilon_1x(t)\big)dt+\alpha_1x(t)d W_1(t),\\
dy(t)=y(t)\big(r_2-\frac{b_2y(t)}{K_2+x(t)}-\varepsilon_2y(t)\big)dt+\alpha_2y(t)d W_2(t).\\
\end{array} \right.
\label{a3}
\end{eqnarray}

On the other hand, population systems may suffer abrupt environmental perturbations, such as epidemics, earthquakes, hurricanes, etc. As a consequence, these systems are very complex and their sample paths may not be continuous, which yields the system \eqref{a3} fail to cope with them. It is recognized that introducing L$\acute{e}$vy noise into the underlying population system may be quite suitable to describe such discontinuous systems. There exists some interesting literatures concerned with SDEs with jumps. We here only mention Bao et al \cite{BMYY, BY}, Liu and Wang \cite{LW4},  Liu and Liang \cite{LL}.  Motivated by those studies, in this paper we consider the following non-autonomous system  with jumps:

\begin{eqnarray}
\left\{
\begin{array}{ll}
dx(t)=x(t^-)\big[\big(r_1(t)-\frac{b_1(t)x(t)}{K_1(t)+y(t)}-\varepsilon_1(t)x(t)\big)dt+\alpha_1(t)d W_1(t)\\ \ +\int_{\mathbb{Y}}\gamma_1(t, u)\tilde{N}(dt,du)\big],\\
dy(t)=y(t^-)\big[\big(r_2(t)-\frac{b_2(t)y(t)}{K_2(t)+x(t)}-\varepsilon_2(t)y(t)\big)dt+\alpha_2(t)d W_2(t)\\ \ +\int_{\mathbb{Y}}\gamma_2(t, u)\tilde{N}(dt,du)\big],\\
\end{array} \right.
\label{a4}
\end{eqnarray}
where $x(t^-)$ and $y(t^-)$ are the left limit of $x(t)$ and $y(t)$ respectively, $r_i(t), b_i(t),$ $ K_i(t), \alpha_i(t), i=1, 2$ are all positive, continuous and bounded functions on $[0, +\infty)$.  N is a Poisson random measure with compensator $\tilde{N}$ and characteristic measure $\mu$
on a measurable subset $\mathbb{Y}$ of $(0, +\infty)$ with $\mu(\mathbb{Y})< +\infty$, $\tilde{N}(dt,du)=N(dt,du)-\mu(du)dt$, $\gamma_i: \mathbb{Y}\times \Omega\rightarrow$ $\mathbb{R}$ is bounded and continuous with respect to $\mu$, and is $\mathfrak{B}(\mathbb{Y})\times \mathfrak{F}_t$-measurable, i=1, 2.

In the next section, the global existence and uniqueness of the positive solution to problem
\eqref{a3} are proved by using comparison theorem for stochastic equations.
Sections 3 is devoted to stochastic boundedness.
Section 4 deals with stochastic permanence. Section 5 discusses the persistence in mean and extinction,
sufficient conditions of persistence in mean and extinction are obtained.\\

Throughout this paper, we let $(\Omega, F, {\mathcal\{F_t\}}_{t\geq 0}, P)$ be a complete probability space with a filtration
${\mathcal\{F_t\}}_{t\geq 0}$ satisfying the usual conditions. For convenience, we assume  that
$X(t)=(x(t), y(t))$ and $|X(t)|=\sqrt{x^2(t)+y^2(t)}.$  $1+\gamma_i(t, u)>0, u\in \mathbb{Y}, i=1, 2,$ there exists a constant $k>0$ such that
$$\int_{\mathbb{Y}}ln(1+\gamma_i(t, u))\vee [ln(1+\gamma_i(t, u))]^2\mu (du)< k, $$
$$\beta_i(t)=0.5\alpha_i^2(t)+\int_{\mathbb{Y}}[\gamma_i(t, u)-ln(1+\gamma_i(t, u))]\mu(du), \  \  i=1, 2,$$
$$Q_i(t)=\int_0^t\int_{\mathbb{Y}}ln(1+\gamma_i(s, u))\tilde{N}(ds, du),  \  \ i=1, 2,$$
$$\hat{f}=\inf_{t\geq 0}f(t),   \  \     \breve{f}=\sup_{t\geq 0}f(t).$$

We end this section by recalling three definitions which we will use in the forthcoming sections.

\begin{defi} $\cite{LM}$ If for any $0<\varepsilon < 1$, there is a constant $\delta(\varepsilon)>0$ such that the solution $X(t)$ of \eqref{a4} satisfies
$$\limsup_{t\rightarrow\infty}P\{{|X(t)|}<\delta\}\geq 1-\varepsilon,$$
for any  initial value $(x_0, y_0)>(0, 0)$, then we say the solution $ X(t)$ be stochastically ultimate boundedness.
\end{defi}

\begin{defi} $\cite{LM}$ If for arbitrary $\varepsilon\in (0,1),$ there are two positive constants $\zeta_1:=\zeta_1(\varepsilon)$ and $\zeta_2:=\zeta_2(\varepsilon)$
such that
$$\liminf_{t\rightarrow \infty}P\{x(t)\leq \zeta_1\}\geq 1-\varepsilon, \ \   \liminf_{t\rightarrow \infty}P\{y(t)\leq \zeta_1\}\geq 1-\varepsilon.\ \ $$
$$\liminf_{t\rightarrow \infty}P\{x(t)\geq \zeta_2\}\geq 1-\varepsilon, \ \   \liminf_{t\rightarrow \infty}P\{y(t)\geq \zeta_2\}\geq 1-\varepsilon.\ \ $$
Then solution of problem \eqref{a4} is said to be stochastically permanent.
\end{defi}

\begin{defi}$\cite{CLCJ}$ If  $x(t), y(t)$ satisfy the following condition
$$\lim_{t\rightarrow \infty}\frac{1}{t}\int_0^tx(s)ds>0, \ \  \lim_{t\rightarrow \infty}\frac{1}{t}\int_0^ty(s)ds>0 \ \  a.s.$$
The problem of \eqref{a4} is said to be persistence in mean.
\end{defi}
%\begin{lem} (Chebyshev's inequality) $\cite{MY}$ If $\delta>0, k>0$ and $X\in L^p(\Omega)$ with $E|X|^k<\infty,$
%then,
%$$P\{|X|\geq \delta\}\leq \delta^{-k}E|X|^k.$$
%\end{lem}

%\begin{lem} $\cite {LW1}.$ Suppose that $x(t)\in C(\Omega\times [0, +\infty), R_+).$\\
%(A) If there exist three positive constants $T, \eta$ and $\eta_0$ such that
%$$ ln x(t)\leq \eta t- \eta_0\int_0^t x(s)ds + \sum_{i=1}^2 \sigma_iW_i(t)+Q(t)$$
%for all  $t\geq T$, where both $\sigma_1$ and $\sigma_2$ are constants, then
%$$\limsup_{t\rightarrow +\infty}\frac{\int_0^t x(s)ds}{t}\leq \frac{\eta}{\eta_0} \  a.s.$$
%(B) If there exist three positive constants $T, \eta$ and $\eta_0$ such that
%$$ ln x(t)\geq \eta t- \eta_0\int_0^t x(s)ds + \sum_{i=1}^2 \sigma_iW_i(t)+Q(t)$$
%for all  $t\geq T$, where both $\sigma_1$ and $\sigma_2$ are constants, then
%$$\liminf_{t\rightarrow +\infty}\frac{\int_0^t x(s)ds}{t}\geq \frac{\eta}{\eta_0} \  a.s.$$
%\end{lem}

\section{Existence and uniqueness of the positive solution}

First, we show that there exists a unique local  positive solution of  \eqref{a4}.
\begin{lem} For the given positive initial value $(x_0, y_0)$, there is  $\tau > 0$ such that
 problem \eqref{a4}
admits a unique positive local solution $X(t)$  a.s. for $t\in [0, \tau).$
\end{lem}
\bpf
  We first set a change of variables : $u(t)=\ln x(t), v(t)=\ln y(t)$, then problem \eqref{a4} deduces to
\begin{eqnarray}
\left\{
\begin{array}{lll}
du(t)=\big(r_1(t)-\beta_1(t)-\frac{b_1(t)e^{u(t)}}{K_1(t)+e^{v(t)}}-\varepsilon_1(t) e^{u(t)}\big)dt+\alpha_1(t) dW_1(t)\\ \ +\int_{\mathbb{Y}}ln(1+\gamma_1(u))\tilde{N}(dt,du),\\
dv(t)=\big(r_2(t)-\beta_2(t)-\frac{b_2(t)e^{v(t)}}{K_2(t)+e^{u(t)}}-\varepsilon_2(t) e^{v(t)}\big)dt+\alpha_2(t) dW_2(t)\\  \ +\int_{\mathbb{Y}}ln(1+\gamma_2(u))\tilde{N}(dt,du)\\
\end{array} \right.
\label{Hb}
\end{eqnarray}
on $t\geq 0$ with initial value $u(0)=\ln{x_0}, v(0)=\ln{y_0} $.
Obviously, the coefficients of \eqref{Hb} satisfy the local Lipschitz condition,
then making use of the theorem \cite{FA, MX} about existence and uniqueness for stochastic differential equation
there is a unique local solution $(u(t), v(t))$ on $t \in [0,\tau )$, where $\tau$ is the explosion time. Hence,
by It$\hat{o}$'s formula, $(x(t), y(t))$ is a unique positive local solution to problem \eqref{a4} with positive initial value.

Next we need to prove solution is  global, that is $\tau=\infty$.
\begin{thm} For any positive initial value $(x_0, y_0)$, there exists a unique global positive solution $(x(t), y(t))$ to problem \eqref{a4}, which
satisfies
$$\lambda (t)\leq x(t)\leq \Lambda(t), \  \theta (t)\leq y(t)\leq \Theta(t), \ t\geq 0, \  a.s.$$
where $\lambda(t),\  \Lambda(t),\  \theta(t)$ and $\Theta(t)$ are defined as \eqref{12}, \eqref{11}, \eqref{14} and \eqref{13}.
\end{thm}
\bpf The reference of \cite{JJS} was the main source of inspiration for its proof. Because of $(x(t), y(t))$ is positive, from the first equation of \eqref{a4}, we can define the following problem
\begin{eqnarray}
\left\{
\begin{array}{l}
d\Lambda(t)=\Lambda(t^-)[\big(r_1(t)-\varepsilon_1(t)\Lambda(t)\big)dt+\alpha_1(t) dW_1(t)+\int_{\mathbb{Y}}\gamma_1(t, u)\tilde{N}(dt,du)],\\
\Lambda(0)=x_0,
\end{array} \right.
\label{HC}
\end{eqnarray}
then
$$
\Lambda(t)=
\frac{e^{\int_0^t(r_1(s)-\beta_1(s))ds+\int_0^t\alpha_1(s)dW_1(s)+Q_1(t)}}{\frac{1}{x_0}+\int^t_0e^{\int_0^s(r_1(u)-\beta_1(u))du+\int_0^s\alpha_1(u)dW_1(u)+Q_1(s)}\varepsilon_1(s) ds}
$$
is the unique solution of \eqref{HC}, and it follows from the comparison theorem for stochastic equations that
\begin{eqnarray}x(t)\leq \Lambda(t),  \  t\in [0,\tau), \ a.s.\label{11}
\end{eqnarray}
On the other hand,
$$
\lambda(t)=\frac{e^{\int_0^t(r_1(s)-\beta_1(s))ds+\int_0^t\alpha_1(s)dW_1(s)+Q_1(t)}}{\frac{1}{x_0}+\int^t_0 e^{\int_0^s(r_1(u)-\beta_1(u))du+\int_0^s\alpha_1(u)dW_1(u)+Q_1(s)}(\varepsilon_1(s)+\frac{b_1(s)}{K_1(s)})ds}
$$
is the solution to the problem
\begin{eqnarray}
\left\{
\begin{array}{llll}
d\lambda(t)=\lambda(t^-)[\big(r_1(t)-(\frac{b_1(t)}{K_1(t)}+\varepsilon_1(t))\lambda(t)\big)dt+\alpha_1(t) dW_1(t)\\ \ +\int_{\mathbb{Y}}\gamma_1(t, u)\tilde{N}(dt,du)],\\
\lambda(0)=x_0,\\
\end{array} \right.
\label{12}
\end{eqnarray}
then
\begin{eqnarray}x(t)\geq \lambda(t),  \  t\in [0,\tau), \ a.s.\label{Hd}
\end{eqnarray}
Similarly, we can get
\begin{eqnarray}y(t)\leq \Theta(t),  \  t\in [0,\tau), \  a.s,\label{13}
\end{eqnarray}
where
$$
\Theta(t)=\frac{e^{\int_0^t(r_2(s)-\beta_2(s))ds+\int_0^t\alpha_2(s)dW_2(s)+Q_2(t)}}{\frac{1}{y_0}+\int^t_0 e^{\int_0^s(r_2(u)-\beta_2(u))du+\int_0^s\alpha_2(u)dW_2(u)+Q_2(s)}\varepsilon_2(s)ds}
$$
and,
\begin{eqnarray}y(t)\geq \theta(t),  \  t\in [0,\tau), \ a. s.\label{14}
\end{eqnarray}
where
$$
\theta(t)=\frac{e^{\int_0^t(r_2(s)-\beta_2(s))ds+\int_0^t\alpha_2(s)dW_2(s)+Q_2(t)}}{\frac{1}{y_0}+\int^t_0 e^{\int_0^s(r_2(u)-\beta_2(u))du+\int_0^s\alpha_2(u)dW_2(u)+Q_2(s)}(\varepsilon_2(s)+\frac{b_2(s)}{K_2(s)})ds}.
$$

Combining \eqref{11}, \eqref{Hd}, \eqref{13} with \eqref{14},  we obtain
$$\lambda (t)\leq x(t)\leq \Lambda(t), \  \theta (t)\leq y(t)\leq \Theta(t), t\geq 0, \  a.s.$$
By Lemma 4.2 in \cite{BMYY}, we know that $\Lambda(t), \lambda(t), \Theta(t), \theta(t)$ will not be exploded in any finite time, it follows from the comparison theorem for stochastic equations \cite{IW}
that $(x(t), y(t))$ exists globally.
\epf

\section{Stochastically ultimate boundedness }
  In a population dynamical system, the nonexplosion property is often not good enough but the property of ultimate boundedness is more desired.
Now, let us present a theorem  about the stochastically ultimate boundedness of \eqref{a4} for any positive initial value.

\begin{thm} Assume that there exists a constant $L(q)>0 $ such that
$$\int_{\mathbb{Y}}|\gamma_i(s, u)|^q\mu(du)\leq \L(q),  \  q>1,  \ i=1, 2.$$
Then  for any positive initial value $(x_0, y_0)$, the solution $X(t)$ of problem \eqref{a4} is stochastically ultimate boundedness.
\end{thm}\
\bpf As the reference of  \cite{BMYY} we define a Lyapunov function $U(x)=x^q$.  By the It$\hat{o}$ formula:
\begin{eqnarray*}
\begin{array}{lll}
E(e^tU(x))&=&U(x_0)+E\int_0^te^s[U(x)+qx^{q-1}dx+\frac{1}{2}q(q-1)x^{q-2}(dx)^2]ds\\
&=&U(x_0)+E\int_0^te^s\{U(x)+q[r_1(s)-\frac{b_1(s)x}{K_1(s)+y}-\varepsilon_1(s)x-\frac{(1-q)\alpha_1^2(s)}{2}\\
&=&+\int_{\mathbb{Y}}[(1+\gamma_1(s, u))^q-1-q\gamma_1(s, u)]\mu(du)]U(x)\}ds\\
&\leq&U(x_0)+E\int_0^te^s\{[-\varepsilon_1(s)x+1+qr_1(s)+\frac{q(q-1)\alpha_1^2(s)}{2}\\
& &+\int_{\mathbb{Y}}[(1+\gamma_1(s, u))^q-1-q\gamma_1(s, u)]\mu(du)]U(x)\}ds.
\end{array}
\end{eqnarray*}
If $ q> 1$, we can deduce that there exists constant $L_1(q)>0$ by assumption such that
\begin{eqnarray*}
\begin{array}{ll}
&U(x)\big\{[(1+qr_1(s)+\frac{q(q-1)}{2}\alpha_1(s)^2)-q\varepsilon_1(t)x]\\
&+\int_{\mathbb{Y}}[(1+\gamma_1(s, u))^q-1-q\gamma_1(s, u)]\mu(du)\big\}\leq L_1(q).\\
\end{array}
\end{eqnarray*}
If $0< q<1$,  using $(1+\gamma_1(s, u))^q-1-q\gamma_1(s, u)\leq 0,$ then we have
\begin{eqnarray*}
\begin{array}{ll}
&U(x)\big\{[(1+qr_1(s)+\frac{q(q-1)}{2}\alpha_1(s)^2)-q\varepsilon_1(t)x]\\
&+\int_{\mathbb{Y}}[(1+\gamma_1(s, u))^q-1-q\gamma_1(s, u)]\mu(du)\big\}\\
&\leq U(x)[1+qr_1(s)-q\varepsilon_1(t)x].\\
\end{array}
\end{eqnarray*}
Therefore,

$$E(e^tU(x))\leq E(U(x_0))+ L_1(q)(e^t-1).$$
Thus,
\begin{eqnarray}
 \limsup_{t\rightarrow\infty} Ex^q\leq L_1(q).
 \label{31}
 \end{eqnarray}
Similarly, we have
\begin{eqnarray}
\limsup_{t\rightarrow\infty}Ey^q\leq L_2(q).
\label{32}
\end{eqnarray}
We now combine \eqref{31}, \eqref{32} with the formula  $\big[x(t)^2+y(t)^2\big]^{\frac{q}{2}}\leq 2^{\frac{q}{2}}\big[x(t)^q+y(t)^q\big]$ to yield
$$\limsup_{t\rightarrow\infty} E|X|^q\leq 2^{\frac{q}{2}}[L_1(q)+L_2(q)]<+\infty.$$
By the Chebyshev's inequality \cite{MX} and the above inequality we can complete the proof.
\epf

\section{Stochastic permanence}

In the study of population models, stochastic permanence is one of the most interesting and important topics. We will discuss this property
by using the method as in \cite{LW4} in this section.
\begin{thm} If  $\min\{\hat{r}_1-\breve{\beta_1}, \hat{r}_2 -\breve{\beta_2}\}>0$, then solution of problem \eqref{a4} is stochastically permanent.
\end{thm}
\bpf For a positive constant $0<\eta <1$, we set a function
$$Z(x)=\frac{1}{x}, V(x)=e^{\lambda t}Z^\eta(x).$$
Straightforward computation  $d V(x)$ by  It$\hat{o}^, s$ formula shows that
\begin{eqnarray*}
\begin{array}{lll}
d V(t)&=&\eta e^{\lambda t}Z^{\eta-2}(t)\big\{-Z^2(t)[r_1(t)-\frac{\alpha_1^2(t)}{2}-\int_{\mathbb{Y}}\gamma_1(t, u)\mu(du)-\frac{(\eta-1)\alpha_1^2(t)}{2}\\
&-&\int_{\mathbb{Y}}\big(\frac{1}{\eta(1+\gamma_1(t, u))^\eta}-\frac{1}{\eta}\big)\mu(du)-\frac{\lambda}{\eta}]+Z^2(t)(\varepsilon_1(t)+\frac{b_1(t)}{K_1(t)+y(t)})\big\}dt\\
&-&\eta \alpha_1 e^{\lambda t}Z^{\eta}(t)dW_1(t)+e^{\lambda t}Z^{\eta}(t)\int_{\mathbb{Y}}\big[(\frac{1}{(1+\gamma_1(t, u))^\eta}-1\big]\tilde{N}(dt, du)\\
&\leq& \eta e^{\lambda t}L Z^{\eta}(t)dt-\eta \alpha_1 e^{\lambda t}Z^{\eta}(t)dW_1(t)+e^{\lambda t}Z^{\eta}(t)\int_{\mathbb{Y}}\big[(\frac{1}{(1+\gamma_1(t, u))^\eta}-1\big]\tilde{N}(dt, du),
\end{array}
\end{eqnarray*}

%\begin{eqnarray*}
%\begin{array}{lll}
%& &d Z(x)=Z(x)[(\alpha_1^2(t)-r_1(t)+\int_{\mathbb{Y}}(\frac{1}{1+\gamma_1(t, u)}-1+\gamma_1(t, u))\mu(du))dt\\
%& &-\alpha_1(t)dW_1(t)+\int_{\mathbb{Y}}(\frac{1}{1+\gamma_1(t, u)}-1)\tilde{N}(dt, du)]+(\varepsilon_1(t)+\frac{b_1(t)}{K_1(t)+y(t)})dt.\\
%\end{array}
%\end{eqnarray*}
Due to
$$-\int_{\mathbb{Y}}ln(1+\gamma_1(t, u))\mu (du)=\lim_{\eta\rightarrow 0^+}\big\{\frac{(\eta-1)\alpha_1^2(t)}{2}+\int_{\mathbb{Y}}\frac{1-(1+\gamma_1(t, u))^\eta}{\eta(1+\gamma_1(t, u))^\eta}\mu(du)\big\},$$
then when $\hat{r}_1 - \breve{\beta_1}>0,$ we can choose a sufficiently small $\eta$ to satisfy
$$r_1(t)-\frac{\alpha_1^2(t)}{2}-\int_{\mathbb{Y}}\gamma_1(t, u)\mu(du)-\big\{\frac{(\eta-1)\alpha_1^2(t)}{2}+\int_{\mathbb{Y}}\frac{1-(1+\gamma_1(t, u))^\eta}{\eta(1+\gamma_1(t, u))^\eta}\mu(du)\big\}>0.$$
Let us choose $\lambda>0$ sufficiently small to satisfy
$$\frac{\lambda}{\eta} < r_1(t)-\frac{\alpha_1^2(t)}{2}-\int_{\mathbb{Y}}\gamma_1(t, u)\mu(du)-\big\{\frac{(\eta-1)\alpha_1^2(t)}{2}+\int_{\mathbb{Y}}\frac{1-(1+\gamma_1(t, u))^\eta}{\eta(1+\gamma_1(t, u))^\eta}\mu(du)\big\}.$$
Then, there is a positive constant $L_1$  satisfying
$$r_1(t)-\frac{\alpha_1^2(t)}{2}-\int_{\mathbb{Y}}\big(\gamma_1(t, u)+\frac{1-(1+\gamma_1(t, u))^\eta}{\eta(1+\gamma_1(t, u))^\eta}\big)\mu(du)$$
$$-\frac{(\eta-1)\alpha_1^2(t)}{2}-\frac{\lambda}{\eta}>-L_1.$$

%Set $V(x)=e^{\lambda t}Z^\eta(x).$ Applying the It$\hat{o}^,s$ formula, we obtain
\begin{eqnarray*}
\begin{array}{lll}
d V(t)
&\leq& \eta e^{\lambda t}L Z^{\eta}(t)dt-\eta \alpha_1 e^{\lambda t}Z^{\eta}(t)dW_1(t)\\
&+& e^{\lambda t}Z^{\eta}(t)\int_{\mathbb{Y}}\big[(\frac{1}{(1+\gamma_1(t, u))^\eta}-1\big]\tilde{N}(dt, du),
\end{array}
\end{eqnarray*}
where $L:=L_1+\breve{\varepsilon_1}+\frac{\breve{b_1}}{\hat{K_1}}.$
Integrating and then taking expectations yields
$$E[V(t)]=e^{\lambda t}E(Z^\eta(x))\leq (\frac{1}{x_0})^\eta+\frac{\eta L}{\lambda}(e^{\lambda t}-1).$$
Therefore,
$$\limsup_{t\rightarrow+\infty}E[\frac{1}{x^\eta(t)}]\leq \frac{\eta L}{\lambda}.$$
Similarly, when $\hat{r}_2 -\breve{\beta_2}>0,$ we have
$$\limsup_{t\rightarrow+\infty}E[\frac{1}{y^\eta(t)}]\leq \frac{\eta L}{\lambda}.$$
For arbitrary $\varepsilon\in(0,1)$, choosing $\zeta_2(\varepsilon)=(\frac{\lambda\varepsilon}{\eta L})^{\frac{1}{\eta}}$  and using Chebyshev inequality, we yield the following inequalities,
$$P\{x(t)<\zeta_2\}=P\{\frac{1}{x^\eta(t)}>\frac{1}{\zeta_2^\eta}\}\leq \frac{E[\frac{1}{x^\eta(t)}]}{\zeta_2^{-\eta}},$$
$$P\{y(t)<\zeta_2\}=P\{\frac{1}{y^\eta(t)}>\frac{1}{\zeta_2^\eta}\}\leq \frac{E[\frac{1}{y^\eta(t)}]}{\zeta_2^{-\eta}}.$$
Hence,
$$\limsup_{t\rightarrow+\infty}P\{x(t)<\zeta_2\}\leq \varepsilon, \ \  \limsup_{t\rightarrow+\infty}P\{y(t)<\zeta_2\}\leq \varepsilon.$$
then,
$$\liminf_{t\rightarrow+\infty}P\{x(t)\geq \zeta_2\}\geq 1-\varepsilon,  \  \   \liminf_{t\rightarrow+\infty}P\{y(t)\geq \zeta_2\}\geq 1-\varepsilon.$$
Combining Chebyshev's inequality with \eqref{31}, \eqref{32}, we can prove that
for arbitrary $\varepsilon\in(0,1)$, there is a positive constant $\zeta_1$ such that
$$\liminf_{t\rightarrow+\infty}P\{x(t)\leq \zeta_1\}\geq 1-\varepsilon,  \ \  \liminf_{t\rightarrow+\infty}P\{x(t)\leq \zeta_1\}\geq 1-\varepsilon.$$
This completes the proof.
\epf

\section{Persistence in mean and extinction}
In the description of population dynamics, it is critical to discuss the property of persistence in mean and extinction.
First, we give a Lemma using the  argument as in \cite{LW1, LW2}
 with suitable modifications.

\begin{lem}  Suppose that $x(t)\in C(\Omega\times [0, +\infty), R_+).$\\
(A) If there exist three positive constants $T, \eta$ and $\eta_0$ such that
$$ ln x(t)\leq \eta t- \eta_0\int_0^t x(s)ds + \int_0^t\sigma_i(s)dW_i(s)+Q(t) \   \   i=1 \ or \ 2$$
for all  $t\geq T$, then
$$\limsup_{t\rightarrow +\infty}\frac{\int_0^t x(s)ds}{t}\leq \frac{\eta}{\eta_0} \  a.s.$$
(B) If there exist three positive constants $T, \eta$ and $\eta_0$ such that
$$ ln x(t)\geq \eta t- \eta_0\int_0^t x(s)ds +\int_0^t\sigma_i(s)dW_i(s)+Q(t) \  \   i=1 \ or \ 2$$
for all  $t\geq T$,  then
$$\liminf_{t\rightarrow +\infty}\frac{\int_0^t x(s)ds}{t}\geq \frac{\eta}{\eta_0} \  a.s.$$
\end{lem}
\bpf (A)
Denote $M_i(t)=\int_0^t\alpha_i(s)dW_i(s)$, $Q_i(t)=\int_0^t\int_{\mathbb{Y}}ln(1+\gamma_i(s, u))\tilde{N}(ds, du),$ then $M_i(t)$, $Q_i(t), i=1, 2$ are real valued local martingales vanishing at $t=0$. One can see that
 the quadratic variations of $M_1(t)$ and $Q_1(t)$ are
 $$\langle M_i(t), M_i(t)\rangle=\int_0^t\alpha_i^2(s)ds\leq \breve{\alpha_i}^2t,$$
 $$\langle Q_i(t), Q_i(t)\rangle=\int_0^t\int_{\mathbb{Y}}(ln(1+\gamma_i(s, u)))^2\mu(du)ds\leq kt, \  i=1, 2 $$
 where $\langle M, M\rangle$ is Meyer's angle bracket process, and
 $$\rho_M(t)=\int_0^t\frac{d\langle M, M\rangle(s)}{(1+s)^2} < max\{k,\breve{\alpha_1}^2\}\int_0^t\frac{ds}{(1+s)^2}<\infty.$$
 By the strong law of large numbers for local martingales \cite{LRA}, we have
 $$\lim_{t\rightarrow\infty}\frac{\int_0^t\alpha_i(s)dW_i(s)}{t} =0, \  \lim_{t\rightarrow\infty}\frac {Q_i(t)}{t}=0, a.s. \  i=1, 2.$$
 Then for arbitrary $\varepsilon>0$, there exists a $T_1>0$ such that for $t>T_1$
 $$ -\varepsilon t<\int_0^t\alpha_i(s)dW_i(s)+Q_i(t)< \varepsilon t.$$
 Set $g(t)=\int_0^t x(s)ds $ for all $t>T_2$, then we have
 $$ln\frac{dg}{dt}\leq (\eta+\varepsilon) t- \eta_0 g, \  t\geq T=max\{T_1, T_2\}.$$
 That is to say : for $t\geq T, e^{\eta_0 g}\frac{dg}{dt}\leq e^{(\eta+\varepsilon)t},$ integrating this inequality from $T$ to $t$, we can get
 $$g(t)\leq \frac{ln\big(e^{\eta_0 T}+\frac{\eta_0}{\eta+\varepsilon}(e^{(\eta+\varepsilon)t}-e^{(\eta+\varepsilon)T}) \big)}{\eta_0}.$$
Therefore
$$\limsup_{t\rightarrow +\infty}\frac{\int_0^t x(s)ds}{t}\leq \limsup_{t\rightarrow +\infty}\frac{ln(\frac{3\eta_0e^{(\eta+\varepsilon)t}}{\eta+\varepsilon})}{\eta_0 t}=\frac{\eta+\varepsilon}{\eta_0}.$$
Using the arbitrariness of $\varepsilon$ we have the assertion.\\
The proof of $(B)$ is similar to $(A)$. The proof is completed.
\epf

 Using Lemma 5.1, we have following theorem.

\begin{thm} Suppose that $\hat{r_i}>\breve{\beta_i},  (i=1, 2), $ $X(t)$ is the positive solution to \eqref{a4} with positive initial value
$(x_0, y_0)$, then the problem \eqref{a4} is persistent in mean.
\end{thm}
\bpf The method is similar to \cite{LL}. We first deduce
$$\limsup_{t\rightarrow \infty}\frac{\ln x(t)}{t}\leq 0, \ \limsup_{t\rightarrow \infty}\frac{\ln y(t)}{t}\leq 0. \ a.s.$$
Making use of It$\hat{o}$'s formula to $e^t lnx$, we deduce
$$e^t lnx-lnx_0=\int_0^t e^s[lnx(s)+r_1(s)-\beta_1(s)-\varepsilon_1(s)x(s)-\frac{b_1(s)x(s)}{K_1(s)+y(s)}]ds+L_1(t)+L_2(t),$$
where $L_1(t)=\int_0^t e^s\alpha_1(s)dW_1(s), \ L_2(t)=\int_0^te^s\int_\mathbb{Y}ln(1+\gamma_1(s,u))\tilde{N}(ds, du)$ are martingales with
the quadratic forms
$$\langle L_1(t), L_1(t)\rangle=\int_0^te^{2s}\alpha_1^2(s)ds,$$
$$\langle L_2(t), L_2(t) \rangle=\int_0^te^{2s}\int_\mathbb{Y}(ln(1+\gamma_1(s,u)))^2\mu(du)ds \leq k\int_0^t e^{2s}ds.$$
By the exponential martingale inequality \cite{MX}, for any positive constants $k, \gamma, \delta,$ we can get that
$$\mathcal{P}\big\{\sup_{0\leq t\leq \gamma k}\big[L_i(t)-0.5e^{-\gamma k}\langle L_i(t), N_i(t) \rangle\big]>\delta e^{\gamma k}lnk\big\}\leq k^{-\delta},$$
it follows from the Borel-Cantelli lemma that for almost all $\omega \in\Omega,$ there is $k_0(\omega)$ such that for each $k\geq k_0(\omega),$
$$L_i(t)\leq 0.5e^{-\gamma k}\langle L_i(t), L_i(t)\rangle+\delta e^{\gamma k} lnk, \ 0\leq t \leq \gamma k.$$
Hence
\begin{eqnarray*}
\begin{array}{lllll}
e^t lnx-lnx_0
&\leq&\int_0^t e^s[ln x(s)+r_1(s)-\beta_1(s)-\varepsilon_(s)x(s)-\frac{b_1(s)x(s)}{K_1(s)+y(s)}]ds\\
&+&\frac{1}{2}e^{-\gamma k}\int_0^te^{2s}\alpha_1^2(s)d s+\frac{1}{2}ke^{-\gamma k}\int_0^t e^{2s}d s+2\delta e^{\gamma k} lnk\\
&=& \int_0^t e^s[ln x(s)+r_1(s)-\int_\mathbb{Y}(\gamma_1(s,u)-ln(1+\gamma_1(s,u)))\mu (du)\\
&-&\varepsilon(s)x(s)
-\frac{b_1(s)x(s)}{K_1(s)+y(s)}-\frac{1}{2}\alpha_1^2(s)[1-e^{s-\gamma k}]-\frac{1}{2}k[1-e^{s-\gamma k}]]ds\\
&+&2\delta e^{\gamma k} lnk\\
&\leq& \int_0^t e^s[ln x(s)+r_1(s)+\int_\mathbb{Y}(|\gamma_1(s,u)|+|ln(1+\gamma_1(s,u))|)\mu (du)\\
&-&\frac{1}{2}\alpha_1^2(s)[1-e^{s-\gamma k}]-\frac{1}{2}k[1-e^{s-\gamma k}]]ds+2\delta e^{\gamma k} lnk .\\
\end{array}
\end{eqnarray*}
Obviously, for any $0\leq s \leq \gamma k$ and $x>0,$ there is a constant A which is independent of k such that
$$ln x(s)+r_1(s)+\int_\mathbb{Y}(|\gamma_1(s,u)|+|ln(1+\gamma_1(s,u))|)\mu (du)$$
$$-0.5\alpha_1^2(s)[1-e^{s-\gamma k}]-0.5k[1-e^{s-\gamma k}]\leq A.$$
Then for $0\leq t \leq \gamma k, k>k_0(\omega),$  we derive
$$e^t lnx-lnx_0\leq A[e^t-1]+2\delta e^{\gamma k} lnk.$$
That is
$$ln x(t)\leq e^{-t}ln x_0+A[1-e^{-t}]+2e^{-t}\delta e^{\gamma k} lnk.$$
Letting $t\rightarrow \infty,$ we have
$$\limsup_{t\rightarrow \infty}\frac{\ln x(t)}{t}\leq 0.$$
Similarly, we get
$$\limsup_{t\rightarrow \infty}\frac{\ln y(t)}{t}\leq 0.$$

On the other hand, applying It$\hat{o}$'s formula to \eqref{12} we have:
$$dln\lambda(t)=(r_1(t)-\beta_1(t)-(\varepsilon_1(t)+\frac{b_1(t)}{K_1(t)})\lambda(t))dt+\alpha_1(t)dW_1(t)$$
$$+\int_{\mathbb{Y}}ln(1+\gamma_1(u))\tilde{N}(dt, du).$$
That is
$$ln\lambda(t)=lnx(0)+\int_0^t\big(r_1(s)-\beta_1(s)-(\varepsilon_1(s)+\frac{b_1(s)}{K_1(s)})\lambda(s)\big)ds$$
$$+\int_0^t\alpha_1(s)W_1(s)ds+Q_1(t).$$
For $t\geq T$, we have
$$ln\lambda(t)\leq (\breve{r}_1 -\hat{\beta_1}+\varepsilon)t-(\hat{\varepsilon}_1 +\frac{\hat{b_1}}{\breve{K_1}})\int_0^t\lambda(s)ds+\int_0^t\alpha_1(s)W_1(s)ds+Q_1(t),$$
$$ln\lambda(t)\geq (\hat{r}_1 -\breve{\beta_1}-\varepsilon)t-(\breve{\varepsilon}_1 +\frac{\breve{b_1}}{\hat{K_1}})\int_0^t\lambda(s)ds+\int_0^t\alpha_1(s)W_1(s)ds+Q_1(t).$$
Let $\varepsilon$ be sufficiently small such that $\hat{r}_1 -\breve{\beta_1}-\varepsilon>0$, then applying Lemma 5.1 to above two inequalities, we get
$$\frac{\hat{K_1}(\hat{r}_1 -\breve{\beta_1}-\varepsilon)}{\breve{b}_1 +\breve{\varepsilon}_1\hat{K_1}}\leq \lim_{t\rightarrow \infty}\inf\frac{\int_0^t\lambda(s)ds}{t}\leq \lim_{t\rightarrow \infty}\sup\frac{\int_0^t\lambda(s)ds}{t}\leq \frac{\breve{K_1}(\breve{r}_1 -\hat{\beta_1}+\varepsilon)}{\hat{b}_1 +\hat{\varepsilon}_1\breve{K_1}}.$$
Making use of the arbitrariness of $\varepsilon$ we get
$$\lim_{t\rightarrow \infty}\sup\frac{\int_0^t\lambda(s)ds}{t}\geq \frac{\hat{K_1}(\hat{r}_1 -\breve{\beta_1}-\varepsilon)}{\breve{b}_1+\breve{\varepsilon}_1\hat{K_1}}.$$
Then
$$ \lim_{t\rightarrow \infty}\sup\frac{\ln \lambda(t)}{t}\geq 0, \  a.s.$$
Therefore
$$\lim_{t\rightarrow \infty}\sup\frac{\ln x(t)}{t}\geq \lim_{t\rightarrow \infty}\sup\frac{\ln \lambda(t)}{t}\geq 0, a.s.$$
To sum up, we have
$$\lim_{t\rightarrow \infty}\frac{\ln x(t)}{t}=0.$$
Similarly, we yield that
$$\lim_{t\rightarrow \infty}\frac{\ln y(t)}{t}=0.$$
Integrating the first equation of \eqref{Hb} from $0$ to $t$, we yield
\begin{eqnarray*}
\begin{array}{llllll}
& &\int_0^t\frac{b_1(s)x(s)}{K_1(s)+y(s)}ds=-\int_0^tdln x(s)++\int_0^t(r_1(s)-\beta_1(s))ds\\ & &+\int_0^t\alpha_1(s)dW_1(s)+Q_1(t)-\int_0^t\varepsilon_1(s)x(s)ds.
\end{array}
\end{eqnarray*}
Because of $\int_0^tb_1(s)x(s)ds \geq \int_0^t\frac{K_1(s)b_1(s)x(s)}{K_1(s)+y(s)}ds,$ we obtain
\begin{eqnarray*}
\begin{array}{lllllll}
&\frac{1}{t}\int_0^tb_1(s)x(s)ds \geq \frac{\hat{K_1}}{t}\big[-(\ln x(t)-\ln x_0)+\int_0^t(r_1(s)-\beta_1(s))ds\\&+\int_0^t\alpha_1(s)dW_1(s)+Q_1(t)-\int_0^t\varepsilon_1(s)x(s)ds\big],
\end{array}
\end{eqnarray*}
which is
$$\frac{1}{t}\int_0^t(b_1(s)+\hat{K_1}\varepsilon_1(s))x(s)ds
 \geq \frac{\hat{K_1}}{t}\big[-(\ln x(t)-\ln x_0)$$
 $$+(\hat{r}_1 -\breve{\beta_1})t+\int_0^t\alpha_1(s)dW_1(s)+Q_1(t)\big].$$
Since that $\lim_{t\rightarrow\infty}\frac{\int_0^t\alpha_1(s)dW_1(s)}{t}=0, \lim_{t\rightarrow\infty}\frac{Q_1(t)}{t}=0$, and $\lim_{t \rightarrow \infty} \frac{\ln x(t)}{t}=0$, we get
$$\lim_{t \rightarrow \infty} \frac{\int_0^tx(s)ds}{t}\geq \frac{\hat{r}_1-\breve{\beta}_1}{\breve{b_1}+\breve{\varepsilon}_1\hat{K_1}}> 0, \ a. s. $$
Similarly, we yield
$$\lim_{t \rightarrow \infty} \frac{\int_0^ty(s)ds}{t}\geq \frac{\hat{r}_2-\breve{\beta_2}}{\breve{b}_2 +\breve{\varepsilon}_2\hat{K_2}}> 0, \ a. s. $$
This completes the proof.
\epf

\begin{thm} Let $X(t)$ be a positive solution of  \eqref{a4} with positive initial value $X(0)$, then\\
$(A)$ If $\breve{r}_1 <\hat{\beta_1}, \breve{r}_2 <\hat{\beta_2}$, then $x(t), y(t)$ be extinction.\\
$(B)$ If $\hat{r}_1 >\breve{\beta_1}, \breve{r}_2 <\hat{\beta_2}$, then $y(t)$ is extinction, $x(t)$ is  persistent in mean.\\
$(C)$ If $\breve{r}_1 <\hat{\beta_1}, \hat{r}_2>\breve{\beta_2}$, then $x(t)$ is extinction, $y(t)$ is  persistent in mean.
 \end{thm}
\bpf We first prove Case $(A)$ of the theorem.
Making use of It$\hat{o}^,$s formula to $ln x, x\in [0, +\infty)$ yields
$$ln x(t)-ln x(0)\leq \int_0^t(r_1(s)-\beta_1(s))ds+\int_0^t\alpha_1(s)dW_1(s)+Q_1(t).$$
% Denote $M_1(t)=\int_0^t\alpha_1(s)dW_1(s)$, then $M_1(t)$, $Q_1(t)$ are real valued local martingales vanishing at $t=0.$ one can see that
% the quadratic variations of $M_1(t)$ and $Q_1(t)$ are
% $$\langle M_1(t), M_1(t)\rangle=\int_0^t\alpha_1^2(s)ds\leq \breve{\alpha_1}t,$$
% $$\langle Q_1(t), Q_1(t)\rangle=\int_0^t\int_{\mathbb{Y}}(ln(1+\gamma_1(s, u)))^2\mu(du)ds\leq kt,$$
% where $\langle M, M\rangle$ is Meyer's angle bracket process. \\
% $$\rho_M(t)=\int_0^t\frac{d\langle M, M\rangle(s)}{(1+s)^2} < max\{k,\breve{\alpha_1}\}\int_0^t\frac{ds}{(1+s)^2}<\infty.$$
% By the strong law of large numbers for local martingales \cite{LRA}, we have
 Because of $$\lim_{t\rightarrow\infty}\frac{\int_0^t\alpha_1(s)dW_1(s)}{t} =0, \  \lim_{t\rightarrow\infty}\frac {Q_1(t)}{t}=0, \; a.s.$$ and $\breve{r}_1 -\hat{\beta_1} < 0$, we can deduce
$$\lim_{t\rightarrow\infty}x(t)=0,  \;   a.s.$$
Similarly
$$\lim_{t\rightarrow\infty}y(t)=0,  \;   a.s.$$
Case $(B).$  Since that  $\breve{r}_2 < \hat{\beta_2}$, we have $\lim_{t\rightarrow\infty}y(t)=0,  \   a.s.$ Then
$$ln x(t)-ln x(0)\leq(\breve{r}_1 - \hat{\beta_1})t-\int_0^t \hat{\varepsilon}_1x(s)ds-\int_0^t\hat{b}_1\frac{x(s)}{\breve{K_1}}ds+\int_0^t\alpha_1(s)dW_1(s)+Q_1(t),$$
$$ln x(t)-ln x(0)\geq(\hat{r}_1 - \breve{\beta_1})t-\breve{\varepsilon}_1\int_0^t x(s)ds-\breve{b}_1 \int_0^t\frac{x(s)}{\hat{K_1}}ds+\int_0^t\alpha_1(s)dW_1(s)+Q_1(t).$$
Making use of Lemma 5.1, we obtain
$$\frac{\hat{K_1}(\hat{r}_1 -\breve{\beta_1})}{\breve{\varepsilon}_1\hat{K_1}+\breve{b}_1}\leq \lim_{t\rightarrow \infty}\inf\frac{\int_0^t x(s)ds}{t}\leq \lim_{t\rightarrow \infty}\sup\frac{\int_0^t x(s)ds}{t}\leq \frac{\breve{K_1}(\breve{r}_1 -\hat{\beta_1})}{\breve{K_1}\hat{\varepsilon}_1+\hat{b}_1}, a.s.$$
Hence,  we get
$$\lim_{t\rightarrow \infty}\inf\frac{\int_0^t x(s)ds}{t}\geq \frac{\hat{K_1}(\hat{r}_1 -\breve{\beta_1})}{\breve{\varepsilon}_1\hat{K_1}+\breve{b}_1} > 0. $$
Case $(C).$
Similar to the arguments in Case $(A)$ and $(B)$, it is easy to find that:\\
$x(t)$ is extinction, $y(t)$ is  persistent in mean, if $\breve{r}_1 <\hat{\beta_1}, \hat{r}_2 >\breve{\beta_2}$.
\epf

\bibliographystyle{ieeetr}
\bibliography{lim}

\bibliography{lim}

\bibliography{lim}
\bibliographystyle{IEEEtran}

\end{document}